\documentclass[11pt]{article}
\usepackage{amssymb}
\usepackage{amsmath}
\usepackage{booktabs}
\usepackage{graphicx}
\usepackage{theorem}

\def\be{{\beta}}

\def\la{{\lambda}}

\def\bbe{{\text{\boldmath $\beta$}}}

\def\bep{{\text{\boldmath $\varepsilon$}}}

\def\bth{{\text{\boldmath $\theta$}}}

\def\bpsi{{\text{\boldmath $\psi$}}}

\def\lah{{\hat \la}}

\def\bbeh{{\widehat \bbe}}

\def\bthh{{\widehat \bth}}

\def\bbet{{\widetilde \bbe}}

\def\Si{{\Sigma}}

\def\La{{\Lambda}}

\def\bSi{{\text{\boldmath $\Si$}}}

\def\bLa{{\text{\boldmath $\La$}}}

\def\bPsi{{\text{\boldmath $\Psi$}}}

\def\bPsih{{\widehat \bPsi}}

\def\u{{\text{\boldmath $u$}}}
\def\v{{\text{\boldmath $v$}}}

\def\y{{\text{\boldmath $y$}}}
\def\z{{\text{\boldmath $z$}}}
\def\A{{\text{\boldmath $A$}}}
\def\B{{\text{\boldmath $B$}}}
\def\C{{\text{\boldmath $C$}}}
\def\D{{\text{\boldmath $D$}}}

\def\G{{\text{\boldmath $G$}}}
\def\H{{\text{\boldmath $H$}}}
\def\I{{\text{\boldmath $I$}}}

\def\X{{\text{\boldmath $X$}}}

\def\Nc{{\cal N}}

\def\tr{{\rm tr\,}}
\def\diag{{\rm diag\,}}

\def\[{{\text{\boldmath $[$}}}
\def\]{{\text{\boldmath $]$}}}
\def\et{{\it et\, al.}}
\def\zero{{\bf\text{\boldmath $0$}}}
\def\one{{\bf\text{\boldmath $1$}}}

\def\|{{\,|\,}}

\def\/{{\Bigr/\!\!}}

\def\1r{{\rm (1)}}
\def\2r{{\rm (2)}}
\def\3r{{\rm (3)}}
\def\4r{{\rm (4)}}
\def\5r{{\rm (5)}}

\def\non{{\nonumber}}
\newtheorem{thm}{Theorem}
\newtheorem{lem}{Lemma}

\theorembodyfont{\rmfamily}

\allowdisplaybreaks[3]
\makeatletter
\newcommand{\Proofname}{Proof}
\makeatother
\makeatletter
\def\BOXSYMBOL{\RIfM@\bgroup\else$\bgroup\aftergroup$\fi
  \vcenter{\hrule\hbox{\vrule height.85em\kern.6em\vrule}\hrule}\egroup}
\makeatother
\newcommand{\BOX}{%
  \ifmmode\else\leavevmode\unskip\penalty9999\hbox{}\nobreak\hfill\fi
  \quad\hbox{\BOXSYMBOL}}

\def\one{{\bf\text{\boldmath $1$}}}
\def\zero{{\bf\text{\boldmath $0$}}}

\begin{document}
\title{On Measuring the Variability of Small Area Estimators in a Multivariate Fay-Herriot Model}

\author{
Tsubasa Ito\footnote{Graduate School of Economics, University of Tokyo, 7-3-1 Hongo, Bunkyo-ku, Tokyo 113-0033, JAPAN. {E-Mail: tsubasa$\_$ito.0710@gmail.com}}
and
Tatsuya Kubokawa\footnote{Faculty of Economics, University of Tokyo, 7-3-1 Hongo, Bunkyo-ku, Tokyo 113-0033, JAPAN. \newline{E-Mail: tatsuya@e.u-tokyo.ac.jp }} 
}

\date{}
\maketitle
\begin{abstract}
This paper is concerned with the small area estimation in the multivariate Fay-Herriot model where  covariance matrix of random effects are fully unknown.
The covariance matrix is estimated by a Prasad-Rao type consistent estimator, and the empirical best linear unbiased predictor (EBLUP) of a vector of small area characteristics is provided.
When the EBLUP is measured in terms of a mean squared error matrix (MSEM),  a second-order approximation of MSEM of the EBLUP and a second-order unbiased estimator of the  MSEM is derived analytically in closed forms.
The performance is investigated through numerical and empirical studies.

\par\vspace{4mm}

\end{abstract}

\section{Introduction}

Mixed effects models and their model-based estimators have been recognized as a useful method in statistical inference.
In particular, small area estimation is an important application of mixed effects models.
Although direct design-based estimates for small area means have large standard errors because of small sample sizes from small areas, the empirical best linear unbiased predictors (EBLUP) induced from mixed effects models provide reliable estimates by ^^ ^^ borrowing strength" from neighboring areas and by using data of auxiliary variables.
Such a model-based method for small area estimation has been studied extensively and actively from both theoretical and applied aspects, mostly for handling univariate survey data.
For comprehensive reviews of small area estimation, see  Ghosh and Rao (1994), Datta and Ghosh (2012), Pfeffermann (2013) and Rao and Molina (2015).

\medskip
When multivariate data with correlations are observed from small areas for estimating multi-dimensional characteristics, like poverty and unemployment indicators, Fay (1987) suggested a multivariate extension of the univariate Fay-Herriot model, called a multivariate Fay-Herriot model, to produce reliable estimates of median incomes for four-, three- and five-person families.
Fuller and Harter (1987) also considered a multivariate modeling for estimating a finite population mean vector.
Datta, Day and Basawa (1999) provided unified theories in empirical linear unbiased prediction or empirical Bayes estimation in general multivariate mixed linear models. 
Datta, Day and Maiti (1998) suggested a hierarchical Bayesian approach to multivariate small area estimation.
Datta, $\et$ (1999) showed the interesting result that the multivariate modeling produces more efficient predictors than the conventional univariate modeling.
Porter, Wikle and Holan (2015) used the multivariate Fay-Herriot model for modeling spatial data.
Ngaruye, von Rosen and Singull (2016) applied a  multivariate mixed linear model to crop yield estimation in Rwanda.

\medskip
Although Datta, $\et$ (1999) developed the general and unified theories concerning the empirical best linear unbiased predictors (EBLUP) and their uncertainty,  it is definitely more helpful and useful to provide concrete forms with closed expressions for EBLUP, the second-order approximation of the mean squared error matrix (MSEM) and the second-order unbiased estimator of the mean squared error matrix.
Recently, Benavent and Morales (2016) treated the multivariate Fay-Herriot model with the covariance matrix of random effects depending on unknown parameters.
As a structure in the covariance matrix, they considered diagonal, AR(1) and the related structures and employed the residual maximum likelihood (REML) method for estimating the unknown parameters embedded in the covariance matrix.
A second-order approximation and estimation of the MESM were also derived.
For some examples, however, we cannot assume specific structures without prior knowledge or information on covariance matrices.  

\medskip
In this paper, we treat the multivariate Fay-Herriot model where the covariance matrix of random effects is fully unknown.
This situation has been studied by Fay (1987), Fuller and Harter (1987), Datta, $\et$ (1998), and useful in the case that statisticians have little knowledge on structures in correlation.
As a specific estimator of the covariance matrix, we employ Prasad-Rao type estimators with closed forms and use the modified versions which are restricted over the space of nonnegative definite matrices.
The empirical best linear unbiased predictors are provided based on the Prasad-Rao type estimators, and second-order approximation of their mean squared error matrices and their second-order unbiased estimators of the MSEM are derived with closed expressions.
These are multivariate extensions of the results given by Prasad and Rao (1990) and Datta, $\et$ (2005) for the univariate case. 

\medskip
The paper is organized as follows: 
Section \ref{sec:EBLUP} gives the Prasad-Rao type estimators and their nonnegative-definite modifications for the covariance matrix of the random effects, and shows their consistency.
In Sections \ref{sec:MSE} and \ref{sec:MSEest}, the second-order approximation of MSEM of EBLUP and the second-order unbiased estimator of the  MSEM are derived in closed forms.
The performance of EBLUP and the MSEM estimator are investigated in Section \ref{sec:sim}.
This numerical study illustrates that the proposals have good performances for the low-dimensional case.
However, a $k\times k$ covariance matrix has  $k(k+1)/2$ parameters, and we need more data so as to maintain the performances of the proposals for higher-dimensional cases.

\medskip
Finally, it is noted that empirical best linear unbiased predictors for small area means are empirical Bayes estimators and related to the so-called James-Stein estimators.
In this sense, the prediction in the multivariate Fay-Herriot model corresponds to the empirical Bayes estimation of a mean matrix of a multivariate normal distribution, which is related to the estimation of a precision matrix from a theoretical aspect as discussed in Efron and Morris (1976).
In this framework, several types of estimators are suggested for estimation of the precision matrix, and it may be an interesting query whether those estimators provide improvements in the multivariate small area estimation.

\section{Empirical Best Linear Unbiased Prediction}
\label{sec:EBLUP}

In this paper, we assume that area-level data $(\y_1, \X_1), \ldots, (\y_m, \X_m)$ are observed, where $m$ is the number of small areas, $\y_i$ is a $k$-variate vector of direct survey estimates and $\X_i$ is a $k\times s$ matrix of covariates associated with $\y_i$ for the $i$-th area.
Then, the multivariate Fay-Herriot model is described as
\begin{equation}
\y_i = \X_i \bbe + \v_i + \bep_i,
\quad i=1, \ldots, m,
\label{eqn:MFH}
\end{equation}
where $\bbe$ is an $s$-variate vector of unknown regression coefficients, $\v_i$ is a $k$-variate vector of random effects depending on the $i$-th area and $\bep_i$ is a $k$-variate vector of sampling errors.
It is assumed that $\v_i$ and $\bep_i$ are mutually independently distributed as
$$
\v_i \sim \Nc_k (\zero, \bPsi)\quad \text{and}\quad \bep_i\sim\Nc_k(\zero, \D_i),
$$
where $\bPsi$ is a $k\times k$ unknown and nonsingular covariance matrix and $\D_1, \ldots, \D_m$ are $k\times k$ known covariance matrices.
This is a multivariate extension of the so-called Fay-Herriot model suggested by Fay and Herriot (1979).

\medskip
For example, we consider the crop data of Battese, Harter and Fuller (1988), who analyze the data in the nested error regression model.
For the $i$-th county, let $y_{i1}$ and $y_{i2}$ be survey data of average areas of corn and soybean, respectively.
Also let $x_{i1}$ and $x_{i2}$ be satellite data of average areas of corn and soybean, respectively.
In this case, $\y_i$, $\X_i$ and $\bbe$ correspond to
$$
\y_i=(y_{i1}, y_{i2})^\top, \quad \X_i=\begin{pmatrix} 1 & x_{i1} & x_{i2} & 0 & 0 & 0\\ 0&0&0&1 & x_{i1} & x_{i2} \end{pmatrix},\quad
\bbe=(\be_1, \ldots, \be_6)^\top
$$
for $k=2$ and $s=6$.
Battese, $\et$ (1988) applied a univariate nested error regression model for each of $y_{i1}$ and $y_{i2}$, while we can use the multivariate model (\ref{eqn:MFH}) for analyzing both data simultaneously.

\medskip
We now express model (\ref{eqn:MFH}) in a matrix form.
Let $\y=(\y_1^\top, \ldots, \y_m^\top)^\top$, $\X=(\X_1^\top, \ldots, \X_m^\top)^\top$, $\v=(\v_1^\top, \ldots, \v_m^\top)^\top$ and $\bep=(\bep_1^\top, \ldots, \bep_m^\top)^\top$.
Then, model (\ref{eqn:MFH}) is expressed as
\begin{equation}
\y=\X\bbe + \v + \bep,
\label{eqn:MFH1}
\end{equation}
where $\v\sim \Nc_{km}(\zero, \I_m \otimes \bPsi)$ and $\bep\sim\Nc_{km}(\zero, \D)$ for $\D=\text{block\ diag}(\D_1, \ldots, \D_m)$.

 \medskip
For the $a$-th area, we want to predict the quantity $\bth_a=\X_a\bbe+\v_a$, which is the conditional mean $E[\y_a\mid \v_a]$ given $\v_a$.
A reasonable estimator can be derived from the conditional expectation $E[\bth_a\mid \y_a]=\X_a\bbe+E[\v_a \mid \y_a]$.
The conditional distribution of $\v_i$ given $\y_i$ and the marginal distribution of $\y_i$ are
\begin{equation}
\begin{split}
\v_i \mid \y_i \sim & \Nc_k(\v_i^*(\bbe, \bPsi), (\bPsi^{-1}+\D_i^{-1})^{-1}),\\
\y_i \sim& \Nc_k(\X_i\bbe, \bPsi+\D_i),
\end{split}
\quad i=1, \ldots, m,
\label{eqn:post}
\end{equation}
where
\begin{equation}
\v_i^*(\bbe,\bPsi) = \bPsi(\bPsi+\D_i)^{-1}(\y_i - \X_i\bbe)
=\big\{ \I_k - \D_i(\bPsi+\D_i)^{-1}\big\}(\y_i-\X_i\bbe).
\label{eqn:vB}
\end{equation}
Thus, we get the estimator
\begin{align*}
\bth_a^*(\bbe, \bPsi)=&\X_a\bbe+E[\v_a \mid \y_a]=\X_a\bbe + \v_a^*(\bbe,\bPsi)\\
=&
\y_a  - \D_a(\bPsi+\D_a)^{-1}(\y_a-\X_a\bbe),
\end{align*}
which corresponds to the Bayes estimator of $\bth_a$ in the Bayesian framework.

\medskip
When $\bPsi$ is known, the maximum likelihood estimator or generalized least squares estimator of $\bbe$ is
\begin{align}
\bbeh(\bPsi) =& 
\{\X^\top (\I_m\otimes \bPsi + \D)^{-1}\X\}^{-1}\X^\top (\I_m\otimes \bPsi + \D)^{-1}\y\non\\
=&\Big\{ \sum_{i=1}^m \X_i^\top(\bPsi+\D_i)^{-1}\X_i\Big\}^{-1} \sum_{i=1}^m \X_i^\top(\bPsi+\D_i)^{-1}\y_i.
\label{eqn:beh}
\end{align}
Substituting $\bbeh(\bPsi)$ into $\bth^*(\bbe,\bPsi)$ yields the estimator
\begin{equation}
\bthh_a(\bPsi) = \y_a  - \D_a(\bPsi+\D_a)^{-1}\big\{\y_a-\X_a\bbeh(\bPsi)\big\}.
\label{eqn:Bayes}
\end{equation}
Datta, $\et$ (1999) showed that $\bthh_a(\bPsi)$ is the best linear unbiased predictor (BLUP) of $\bth_a$.
It can be also demonstrated that $\bthh_a(\bPsi)$ is the Bayes estimator against the uniform prior distribution of $\bbe$ as well as the empirical Bayes estimator as shown above, which is called the Bayes empirical Bayes estimator.

\medskip
Concerning  estimation of $\bPsi$, it is noted that $E[(\y_i-\X_i\bbe)(\y_i-\X_i\bbe)^\top]=\bPsi+\D_i$ for $i=1, \ldots, m$, which implies that $\sum_{i=1}^m E[(\y_i-\X_i\bbe)(\y_i-\X_i\bbe)^\top]=m \bPsi+\sum_{i=1}^m\D_i$.
Substituting the ordinary least squares estimator $\bbet=(\X^\top\X)^{-1}\X^\top\y$ into $\bbe$, we get the consistent estimator
\begin{equation}
\bPsih_0 = {1\over m} \sum_{i=1}^m \big\{ (\y_i-\X_i\bbet)(\y_i-\X_i\bbet)^\top-\D_i\big\}.
\label{eqn:Psi0}
\end{equation}

Taking the expectation of $\bPsih_0$, we can see that $E[\bPsih_0]=\bPsi+{\rm Bias}_{\bPsih_0}(\bPsi)$, where
\begin{align}
{\rm Bias}_{\bPsih_0}(\bPsi)=&
 {1\over m}\sum_{i=1}^m \X_i(\X^\top\X)^{-1}\Big\{\sum_{j=1}^m\X_j^\top (\bPsi+\D_j)\X_j\Big\}(\X^\top\X)^{-1} \X_i^\top
\non\\
&
-{1\over m}\sum_{i=1}^m (\bPsi+\D_i)\X_i(\X^\top\X)^{-1}\X_i^\top
-{1\over m}\sum_{i=1}^m \X_i(\X^\top\X)^{-1}\X_i^\top(\bPsi+\D_i).
\label{eqn:Bias}
\end{align}
Substituting $\bPsih_0$ into ${\rm Bias}_{\bPsih_0}(\bPsi)$, we get a bias-corrected given by
\begin{equation}
\bPsih_1 = \bPsih_0 - {\rm Bias}_{\bPsih_0}(\bPsih_0).
\label{eqn:Psi1}
\end{equation}

For notational convenience, we use the same notation $\bPsih$ for $\bPsih_0$ and $\bPsih_1$ without any confusion.
It is noted that both estimators are not necessarily nonnegative definite.
In this case, there exist a $k\times k$ orthogonal matrix $\H$ and a diagonal matrix $\bLa=\diag(\la_1, \ldots, \la_k)$ such that $\bPsih=\H\bLa \H^\top$.
Let $\bLa^+=\diag(\max\{0, \la_1\}, \ldots, \max\{0, \la_k\})$, and let
$$
\bPsih^+ = \H\bLa^+\H^\top.
$$
Replace $\bPsi$ in $\bthh_a(\bPsi)$ with the estimator $\bPsih^+$, and the resulting estimator is the empirical Bayes (EB) estimator
\begin{equation}
\bthh_a^{EB} = \bthh_a(\bPsih^+).
\label{eqn:EB}
\end{equation}

To guarantee asymptotic properties of $\bPsih^+$, we assume the following conditions:

(H1)\ $0<k<\infty$, $0<s<\infty$.

(H2)\ There exist positive constants ${\underline d}$ and ${\overline d}$ such that ${\underline d}$ and ${\overline d}$ do not depend on $m$ and satify ${\underline d}\I_k \leq \D_i \leq {\overline d}\I_k$ for $i=1, \ldots, m$.

(H3)\ $\X^\top \X$ is nonsingular and $\X^\top \X/m$ converges to a  positive definite matrix.

\begin{thm}
\label{thm:BLUP}
Under conditions {\rm (H1)}-{\rm (H3)}, the following properties hold for $\bPsih=\bPsih_0$ and $\bPsih_1$:

{\rm (1)}\ ${\rm Bias}_{\bPsih_0}(\bPsi)=O(m^{-1})$, which means that $\bPsih_0$ has the second-order bias, while $\bPsih_1$ is a second-order unbiased estimator of $\bPsi$.

{\rm (2)}\ $\bPsih-\bPsi=O_p(m^{-1/2})$ and $\bbeh(\bPsih)-\bbe=O_p(m^{-1/2})$.

{\rm (3)}\ The nonnegative defnite matrix $\bPsih^+$ is consistent for large $m$, and $P(\bPsih^+\not= \bPsih)=O(m^{-K})$ for any $K$.
\end{thm}

{\bf Proof}.\ \ We begin with writing $\bPsih_0-\bPsi$ as
\begin{align*}
\bPsih_0-\bPsi=&
{1\over m}\sum_{i=1}^m\{(\y_i-\X_i\bbe)(\y_i-\X_i\bbe)^\top - (\bPsi+\D_i)\}
+ {1\over m}\sum_{i=1}^m \X_i(\bbet-\bbe)(\bbet-\bbe)^\top\X_i^\top
\non\\
&
-{1\over m}\sum_{i=1}^m (\y_i-\X_i\bbe)(\bbet-\bbe)^\top\X_i^\top
-{1\over m}\sum_{i=1}^m \X_i(\bbet-\bbe)(\y_i-\X_i\bbe)^\top,
\end{align*}
which yields the bias given in (\ref{eqn:Bias}).
It is easy to check that the bias is of order $O(m^{-1})$.

\medskip
For (2), it is noted that $\bPsih-\bPsi$ is approximated as
\begin{align}
\bPsih-\bPsi=&
{1\over m}\sum_{i=1}^m\{(\y_i-\X_i\bbe)(\y_i-\X_i\bbe)^\top - (\bPsi+\D_i)\}
+O_p(m^{-1})
\non\\
=& {1\over m}\sum_{i=1}^m\{\u_i\u_i^\top - (\bPsi+\D_i)\}+O_p(m^{-1}),
\label{eqn:Psiap}
\end{align}
where $\u_i=\y_i-\X_i\bbe$, having  $\Nc_k(\zero, \bPsi+\D_i)$.
It is here noted that $(\u_i\u_i^\top - (\bPsi+\D_i))/m$ for $i=1,\ldots,m$ are mutually  independent and $E(\u_i\u_i^\top - (\bPsi+\D_i))/m=0$ for $i=1,\ldots,m$.
Then the consistency follows because $\sum_{i=1}^mE(\u_i\u_i^\top - (\bPsi+\D_i))^2/m^2=\sum_{i=1}^m(2(\bPsi+\D_i)^2+\tr(\bPsi+\D_i)\I_{k})/m^2=O(m^{-1})$ under condition (H2).
Using condition (H2) and finiteness of moments of normal random variables, we can show that $\sqrt{m}(\bPsih-\bPsi)$ converges to a multivariate normal distribution, which implies that $\bPsih-\bPsi=O_p(m^{-1/2})$.

\medskip
We next verify that $\bbeh(\bPsih)-\bbe=O_p(m^{-1/2})$.
Note that $\bbeh(\bPsih)-\bbe$ is decomposed as $\{\bbeh(\bPsih)-\bbeh(\bPsi)\}+\{\bbeh(\bPsi)-\bbe\}$.
For $\bbeh(\bPsi)-\bbe$, it is noted that
\begin{align}
\bbeh(\bPsi)-\bbe=
\Big\{ \sum_{i=1}^m \X_i^\top(\bPsi+\D_i)^{-1}\X_i\Big\}^{-1} \sum_{i=1}^m \X_i^\top(\bPsi+\D_i)^{-1}(\y_i-E\y_i).
\end{align}
Then, ${\rm Var}(\bbeh(\bPsi)-\bbe)=\Big\{ \sum_{i=1}^m \X_i^\top(\bPsi+\D_i)^{-1}\X_i\Big\}^{-1}=O(1/m)$ and this implies $\bbeh(\bPsi)-\bbe=O_{p}(m^{-1/2})$.
We next evaluate $\bbeh(\bPsih)-\bbe(\bPsi)$ as
\begin{align}
&\bbeh(\bPsih)-\bbe(\bPsi)\non\\
&=\Big\{ \sum_{i=1}^m \X_i^\top(\bPsih+\D_i)^{-1}\X_i\Big\}^{-1} \sum_{i=1}^m \X_i^\top(\bPsih+\D_i)^{-1}\y_i
\non\\
&
-\Big\{ \sum_{i=1}^m \X_i^\top(\bPsi+\D_i)^{-1}\X_i\Big\}^{-1} \sum_{i=1}^m \X_i^\top(\bPsi+\D_i)^{-1}\y_i
\non\\
&=
\Big\{ \sum_{i=1}^m \X_i^\top(\bPsih+\D_i)^{-1}\X_i\Big\}^{-1} \sum_{i=1}^m \X_i^\top\Big\{(\bPsih+\D_i)^{-1}-(\bPsi+\D_i)^{-1}\Big\}\y_i
\non\\
&
+\Big[ \Big\{ \sum_{i=1}^m \X_i^\top(\bPsih+\D_i)^{-1}\X_i\Big\}^{-1}-\Big\{ \sum_{i=1}^m \X_i^\top(\bPsi+\D_i)^{-1}\X_i\Big\}^{-1} \Big]
 \sum_{i=1}^m \X_i^\top(\bPsi+\D_i)^{-1}\y_i
\non\\
&=
I_{1}+I_{2}.
\label{eqn:pbe1}
\end{align}
First, $I_{1}$ is written as
\begin{equation}
I_{1}=
-\Big\{ \sum_{i=1}^m \X_i^\top(\bPsih+\D_i)^{-1}\X_i\Big\}^{-1} \sum_{i=1}^m \X_i^\top(\bPsih+\D_i)^{-1}(\bPsih-\bPsi)(\bPsi+\D_i)^{-1}\y_i,
\label{eqn:pbe2}
\end{equation}
which is of order $O_p(m^{-1/2})$, because $\sum_{i=1}^m \X_i^\top(\bPsih+\D_i)^{-1}\X_i=O_{p}(m)$ and $\sum_{i=1}^m \X_i^\top(\bPsih+\D_i)^{-1}(\bPsih-\bPsi)(\bPsi+\D_i)^{-1}\y_i=O_{p}(m^{1/2})$.
Next, $I_{2}$ is rewritten as
\begin{align}
I_{2}=&
-\Big\{ \sum_{i=1}^m \X_i^\top(\bPsih+\D_i)^{-1}\X_i\Big\}^{-1} \sum_{i=1}^m \X_i^\top \Big\{ (\bPsih+\D_i)^{-1}-(\bPsi+\D_i)^{-1} \Big\} \X_i
\non\\
&
\times \Big\{ \sum_{i=1}^m \X_i^\top(\bPsi+\D_i)^{-1}\X_i\Big\}^{-1}\sum_{i=1}^m \X_i^\top(\bPsi+\D_i)^{-1}\y_i
\non\\
=&
\Big\{ \sum_{i=1}^m \X_i^\top(\bPsih+\D_i)^{-1}\X_i\Big\}^{-1} \sum_{i=1}^m \X_i^\top(\bPsih+\D_i)^{-1}(\bPsih-\bPsi)(\bPsi+\D_i)^{-1}\X_i
\non\\
&
\times \Big\{ \sum_{i=1}^m \X_i^\top(\bPsi+\D_i)^{-1}\X_i\Big\}^{-1}\sum_{i=1}^m \X_i^\top(\bPsi+\D_i)^{-1}\y_i,
\label{eqn:pbe3}
\end{align}
which is of order $O_p(m^{-1/2})$, because $\sum_{i=1}^m \X_i^\top(\bPsih+\D_i)^{-1}\X_i=O_{p}(m)$, $\sum_{i=1}^m \X_i^\top(\bPsih+\D_i)^{-1}(\bPsih-\bPsi)(\bPsi+\D_i)^{-1}\X_i=O_{p}(m^{1/2})$, $\sum_{i=1}^m \X_i^\top(\bPsi+\D_i)^{-1}\X_i=O(m)$ and $\sum_{i=1}^m \X_i^\top(\bPsih+\D_i)^{-1}\y_i=O_{p}(m)$. 
Thus, we have $\bbeh(\bPsih)-\bbe(\bPsi)=O_p(m^{-1/2})$, and it is concluded that $\bbeh(\bPsih)-\bbe=O_{p}(m^{-1/2})$.

\medskip
For (3), let $\lah_1, \ldots, \lah_k$ be eigenvalues of $\bPsih$, and let $\la_1, \ldots, \la_k$ be eigenvalues of $\bPsi$.
Then, for $j=1, \ldots, k$, 
$$
P(\lah_j < 0)=P(\lah_j-\la_j < - \la_j)=P( -(\lah_j-\la_j) > \la_j)
\leq P(|\sqrt{m}(\lah_j-\la_j)| > \sqrt{m}\la_j).
$$
Note that $\la_j>0$.
It follows from the Markov inequality that for any $K>0$, 
$$
P(|\sqrt{m}(\lah_j-\la_j)| > \sqrt{m}\la_j)
\leq {E[\{|\sqrt{m}(\lah_j-\la_j)|\}^{2K}] \over  (\sqrt{m}\la_j)^{2K}}=O(m^{-K}),
$$
because $\lah_j-\la_j=O_p(m^{-1/2})$ from $\bPsih-\bPsi=O_p(m^{-1/2})$.
\hfill$\Box$

\section{Second-order Approximation of Mean Squared Error Matrix}
\label{sec:MSE}

Uncertainty of the empirical Bayes estimator $\bthh_a^{EB}$ in (\ref{eqn:EB}) is measured by the mean squared error matrix (MSEM), defined as ${\rm MSEM}(\bthh_a^{EB})=E[\{\bthh_a^{EB}-\bth_a\}\{\bthh_a^{EB}-\bth_a\}^\top]$.
It is noted that
$$
\bthh_a^{EB}-\bth_a
=
\{\bth_a^*(\bbe,\bPsi)-\bth_a\} + \{\bthh_a(\bPsi)-\bth_a^*(\bbe,\bPsi)\}
+\{\bthh_a^{EB}-\bthh_a(\bPsi)\}
$$
and that $\bthh_a(\bPsi)-\bth_a^*(\bbe,\bPsi)=-D_a(\bPsi+D_a)^{-1}\X_a\{\bbeh(\bPsi)-\bbe\}$.
The following lemma is useful for evaluating the mean square error matrix.

\begin{lem}
\label{lem:1}
$\bbeh(\bPsi)$ is independent of $\y-\X\bbet$ or $\bPsih$.
Also, $\bbeh(\bPsi)$ is independent of $\bthh_a^{EB}-\bthh_a(\bPsi)$.
\end{lem}

{\bf Proof}.\ \ 
The covariance of $\y-\X\bbet$ and $\bbeh(\bPsi)$ is
\begin{align*}
E[&(\y-\X\bbet)(\bbeh(\bPsi)-\bbe)^\top] \{\X^\top (\I_m\otimes \bPsi+\D)^{-1}\X\} \\
&=E[(\y-\X\bbet)(\y-\X\bbe)^\top](\I_m\otimes \bPsi+\D)^{-1}\X
\\
&=\Big[(\I_m\otimes \bPsi+\D) - \X \{\X^\top (\I_m\otimes \bPsi+\D)^{-1}\X\}^{-1}\X^\top \Big](\I_m\otimes \bPsi+\D)^{-1}\X
\\
&=\zero.
\end{align*}
This implies that $\bbeh(\bPsi)$ is independent of $\y-\X\bbet$ or $\bPsih$.
It is also noted that
\begin{align*}
\bthh_a^{EB}-\bthh_a(\bPsi)
=&
-\D_a(\bPsih+\D_a)^{-1}(\y_a-\X_a\bbeh(\bPsih))+\D_a(\bPsi+\D_a)^{-1}(\y_a-\X_a\bbeh(\bPsi))
\\
=&
\D_a\{(\bPsi+\D_a)^{-1}-(\bPsih+\D_a)^{-1}\}(\y_a-\X_a\bbet)
\\
&+\D_a(\bPsih+\D_a)^{-1}\X_a\{\X^\top (\I_m\otimes \bPsih+\D)^{-1}\X\}^{-1}\X^\top (\I_m\otimes \bPsih+\D)^{-1}(\y-\X\bbet)
\\
&- \D_a(\bPsi+\D_a)^{-1}\X_a\{\X^\top (\I_m\otimes \bPsi+\D)^{-1}\X\}^{-1}\X^\top (\I_m\otimes \bPsi+\D)^{-1}(\y-\X\bbet),
\end{align*}
which is a function of $\y-\X\bbet$.
Hence, $\bthh_a^{EB}-\bthh_a(\bPsi)$ is independent of $\bbeh(\bPsi)$.
\hfill$\Box$

\bigskip 
Using Lemma \ref{lem:1}, we can decompose the mean squared error matrix as
\begin{align}
{\rm MSEM}(\bthh_a^{EB})=&
E[\{\bth_a^*(\bbe,\bPsi)-\bth_a\}\{\bth_a^*(\bbe,\bPsi)-\bth_a\}^\top]
\non\\
&+E[\{\bthh_a(\bPsi)-\bth_a^*(\bbe,\bPsi)\}\{\bthh_a(\bPsi)-\bth_a^*(\bbe,\bPsi)\}^\top]
\non\\
&+E[\{\bthh_a^{EB}-\bthh_a(\bPsi)\}\{\bthh_a^{EB}-\bthh_a(\bPsi)\}^\top]
\non\\
=& \G_{1a}(\bPsi) + \G_{2i}(\bPsi) + E[\{\bthh_a^{EB}-\bthh_a(\bPsi)\}\{\bthh_a^{EB}-\bthh_a(\bPsi)\}^\top],
\label{eqn:MSE}
\end{align}
where
\begin{equation}
\begin{split}
\G_{1a}(\bPsi)=&
(\bPsi^{-1}+\D_a^{-1})^{-1}=\bPsi(\bPsi+\D_a)^{-1}\D_a,\\
\G_{2a}(\bPsi)=&
\D_a(\bPsi+\D_a)^{-1}\X_a\{\X^\top (\I_m\otimes \bPsi+\D)^{-1}\X\}^{-1}\X_a^\top(\bPsi+\D_a)^{-1}\D_a.
\end{split}
\label{eqn:G12}
\end{equation}
The third term can be approximated as
\begin{align}
\G_{3a}(\bPsi)=&
{1\over m^2}\D_a(\bPsi+\D_a)^{-1}
\Big[ \sum_{i=1}^m(\bPsi+\D_i)(\bPsi+\D_a)^{-1}(\bPsi+\D_i) 
\non\\
&+ \sum_{i=1}^m\{\tr[(\bPsi+\D_i) (\bPsi+\D_a)^{-1}]\}(\bPsi+\D_i)\Big](\bPsi+\D_a)^{-1}\D_a.
\label{eqn:G3}
\end{align}

\begin{thm}
\label{thm:MSE}
The mean squared error matrix of the empirical Bayes estimator $\bthh_a^{EB}$ is approximated as
\begin{equation}
{\rm MSEM}(\bthh_a^{EB})=
\G_{1a}(\bPsi) + \G_{2a}(\bPsi) + \G_{3a}(\bPsi) + O(m^{-3/2}).
\label{eqn:MSE}
\end{equation}
\end{thm}

{\bf Proof}.\ \ 
We shall prove that $E[\{\bthh_a^{EB}-\bthh_a(\bPsi)\}\{\bthh_a^{EB}-\bthh_a(\bPsi)\}^\top]=\G_{3a}(\bPsi)+O_p(m^{-3/2})$.
Also from (2) in Theorem \ref{thm:BLUP}, it is sufficient to show this approximation for $\bPsih$ instead of $\bPsih^+$.
It is observed that
\begin{align*}
\bthh_a^{EB}-\bthh_a(\bPsi)
=& \D_a\{(\bPsi+\D_a)^{-1}-(\bPsih+\D_a)^{-1}\}(\y_a-\X_a\bbe)
+\D_a(\bPsih+\D_a)^{-1}\X_a\{\bbeh(\bPsih)-\bbe\}
\\
&-\D_a(\bPsi+\D_a)^{-1}\X_a\{\bbeh(\bPsi)-\bbe\}.
\end{align*}
Using the equation 
\begin{equation}
(\bPsih+\D_i)^{-1}=(\bPsi+\D_i)^{-1}-(\bPsi+\D_i)^{-1}(\bPsih-\bPsi)(\bPsih+\D_i)^{-1},
\label{eqn:identity}
\end{equation}
we can see that
\begin{align*}
\D_a&\{(\bPsi+\D_a)^{-1}-(\bPsih+\D_a)^{-1}\}(\y_a-\X_a\bbe)\\
=& \D_a(\bPsi+\D_a)^{-1}(\bPsih-\bPsi)(\bPsih+\D_a)^{-1}(\y_a-\X_a\bbe)
\\
=&\D_a(\bPsi+\D_a)^{-1}(\bPsih-\bPsi)(\bPsi+\D_a)^{-1}(\y_a-\X_a\bbe)+O_p(m^{-1})
\end{align*}
and
\begin{align*}
\D_a&(\bPsih+\D_a)^{-1}\X_a\{\bbeh(\bPsih)-\bbe\}\\
=&\D_a(\bPsi+\D_a)^{-1}\X_a\{\bbeh(\bPsih)-\bbe\}
- \D_a(\bPsi+\D_a)^{-1}(\bPsih-\bPsi)(\bPsih+\D_a)^{-1}\X_a\{\bbeh(\bPsih)-\bbe\}
\\
=&\D_a(\bPsi+\D_a)^{-1}\X_a\{\bbeh(\bPsih)-\bbe\}+O_p(m^{-1}).
\end{align*}
Thus, we have
\begin{align*}
\bthh_a^{EB}-\bthh_a(\bPsi)
=&\D_a(\bPsi+\D_a)^{-1}(\bPsih-\bPsi)(\bPsi+\D_a)^{-1}(\y_a-\X_a\bbe)\\
&
+ \D_a(\bPsi+\D_a)^{-1}\X_a\{\bbeh(\bPsih)-\bbeh(\bPsi)\}+O_p(m^{-1})\\
=&I_1+I_2 + O_p(m^{-1}). \quad \text{(say)}
\end{align*}
For $I_2$, it is noted that
\begin{align*}
\bbeh&(\bPsih)-\bbeh(\bPsi)\\
=& \Big[\Big\{\sum_{j=1}^m \X_j^\top(\bPsih+\D_j)^{-1}\X_j\Big\}^{-1}-\Big\{\sum_{j=1}^m \X_j^\top(\bPsi+\D_j)^{-1}\X_j\Big\}^{-1}\Big] 
\\
&\times \sum_{i=1}^m \X_i^\top(\bPsih+\D_i)^{-1}(\y_i-\X_i\bbe)
\\
&+\Big\{\sum_{j=1}^m \X_j^\top(\bPsi+\D_j)^{-1}\X_j\Big\}^{-1}\sum_{i=1}^m \X_i^\top\Big\{(\bPsih+\D_i)^{-1}-(\bPsi+\D_i)^{-1}\Big\}(\y_i-\X_i\bbe)\\
=&I_{21}+I_{22}.
\end{align*}
We can evaluate $I_{21}$ as
\begin{align*}
I_{21}=&\Big\{\sum_{j=1}^m \X_j^\top(\bPsi+\D_j)^{-1}\X_j\Big\}^{-1}\Big\{\sum_{i=1}^m \X_i^\top (\bPsi+\D_i)^{-1}(\bPsih-\bPsi)(\bPsi+\D_i)^{-1}\X_i \Big\} \{\bbeh(\bPsih)-\bbe\}
\\
=& O_p(m^{-1}),
\end{align*}
because $\sum_{j=1}^m \X_j^\top(\bPsi+\D_j)^{-1}\X_j=O(m)$, $\sum_{i=1}^m \X_i^\top (\bPsi+\D_i)^{-1}(\bPsih-\bPsi)(\bPsi+\D_i)^{-1}\X_i =O_p(m^{1/2})$ and $\bbeh(\bPsih)-\bbe=O_p(m^{-1/2})$ from Theorem \ref{thm:BLUP} (2).
We next estimate $I_{22}$ as
\begin{align*}
I_{22}=&
- \Big\{\sum_{j=1}^m \X_j^\top(\bPsi+\D_j)^{-1}\X_j\Big\}^{-1}
\Big\{\sum_{i=1}^m \X_i^\top\A(\bPsih,\D_i)\X_i\Big\}\\
&\times \Big\{\sum_{i=1}^m \X_i^\top \A(\bPsih,\D_i)\X_i\Big\}^{-1}\sum_{i=1}^m \X_i^\top\A(\bPsih,\D_i)(\y_i-\X_i\bbe)
\end{align*}
for $\A(\bPsih,\D_i)=(\bPsih+\D_i)^{-1}(\bPsih-\bPsi)(\bPsi+\D_i)^{-1}$.
It can be seen that $I_{22}=O_p(m^{-1})$ from the same arguments as in $I_{21}$.
Thus, it follows that $I_2=O_p(m^{-1})$.
Hence, we have
\begin{align*}
E[\{&\bthh_a^{EB}-\bthh_a(\bPsi)\}\{\bthh_a^{EB}-\bthh_a(\bPsi)\}^\top]
\\
=&\D_a(\bPsi+\D_a)^{-1}E\Big[(\bPsih-\bPsi)(\bPsi+\D_a)^{-1}(\y_a-\X_a\bbe)
(\y_a-\X_a\bbe)^\top (\bPsi+\D_a)^{-1}(\bPsih-\bPsi)\Big]
\\
&\times  (\bPsi+\D_a)^{-1}\D_a + O(m^{-3/2}).
\end{align*}

It is noted from (\ref{eqn:Psiap}) that $\bPsih-\bPsi$ is approximated as
$$
\bPsih-\bPsi = {1\over m}\sum_{i=1}^m\{\u_i\u_i^\top - (\bPsi+\D_i)\}+O_p(m^{-1}),
\label{eqn:Psiap}
$$
which is used to evaluate 
\begin{align*}
E\Big[&(\bPsih-\bPsi)(\bPsi+\D_a)^{-1}(\y_a-\X_a\bbe)
(\y_a-\X_a\bbe)^\top (\bPsi+\D_a)^{-1}(\bPsih-\bPsi)\Big]
\\
=&
{1\over m^2}\sum_{i=1}^m\sum_{j=1}^m 
E\Big[ \{\u_i\u_i^\top-(\bPsi+\D_i)\}(\bPsi+\D_a)^{-1}\u_a\u_a^\top (\bPsi+\D_a)^{-1}
\{\u_j\u_j^\top-(\bPsi+\D_i)\}\Big]\\
&+O(m^{-3/2})\non\\
=&
{1\over m^2}\sum_{i=1}^m
E\Big[ \{\u_i\u_i^\top-(\bPsi+\D_i)\}(\bPsi+\D_a)^{-1}\u_a\u_a^\top (\bPsi+\D_a)^{-1}
\{\u_i\u_i^\top-(\bPsi+\D_i)\}\Big]\\
&+O(m^{-3/2}),
\end{align*}
since $E\Big[ \{\u_i\u_i^\top-(\bPsi+\D_i)\}(\bPsi+\D_a)^{-1}\u_a\u_a^\top (\bPsi+\D_a)^{-1}\{\u_j\u_j^\top-(\bPsi+\D_i)\}\Big]=\zero$ for $i\not= j$.
Letting $\z_i=(\bPsi+\D_i)^{-1/2}\u_i$, we can see that  $\z_i\sim \Nc_k(\zero, \I_k)$.
Then,
\begin{align*}
{1\over m^2}&\sum_{i=1}^m
E\Big[ \{\u_i\u_i^\top-(\bPsi+\D_i)\}(\bPsi+\D_a)^{-1}\u_a\u_a^\top (\bPsi+\D_a)^{-1}
\{\u_i\u_i^\top-(\bPsi+\D_i)\}\Big]\\
=&
{1\over m^2}\sum_{i\not= a}
(\bPsi+\D_i)^{1/2}E\Big[ (\z_i\z_i^\top-\I)\B\z_a\z_a^\top \B^\top (\z_i\z_i^\top-\I)\Big](\bPsi+\D_i)^{1/2} + O(m^{-2}),
\end{align*}
for $\B=(\bPsi+\D_i)^{1/2}(\bPsi+\D_a)^{-1/2}$.
Let $\C=\B\B^\top=(\bPsi+\D_i)^{1/2} (\bPsi+\D_a)^{-1} (\bPsi+\D_i)^{1/2}$.
For $i\not= a$, 
\begin{align*}
E[& (\z_i\z_i^\top-\I)\B\z_a\z_a^\top \B^\top (\z_i\z_i^\top-\I)]
\\
=&E[ \z_i\z_i^\top\B\z_a\z_a^\top \B^\top \z_i\z_i^\top+\B\z_a\z_a^\top \B^\top
-\z_i\z_i^\top\B\z_a\z_a^\top \B^\top -\B\z_a\z_a^\top \B^\top \z_i\z_i^\top]
\\
=&
E[ \z_i\z_i^\top\C \z_i\z_i^\top-\C]
=\C+(\tr\C)\I_k,
\end{align*}
because $E[\z_i\z_i^\top\C \z_i\z_i^\top]=2\C+(\tr\C)\I_k$.
Thus,
\begin{align*}
{1\over m^2}&\sum_{i\not= a}
(\bPsi+\D_i)^{1/2}\{ \C+(\tr\C)\I_k\}(\bPsi+\D_i)^{1/2}
\\
=&
{1\over m^2}\sum_{i=1}^m
(\bPsi+\D_i)^{1/2}\{ \C+(\tr\C)\I_k\}(\bPsi+\D_i)^{1/2} + O(m^{-2}),
\end{align*}
which leads to the expression in (\ref{eqn:G3}).
\hfill$\Box$

\section{Estimation of Mean Squared Error Matrix}
\label{sec:MSEest}

In this section, we obtain a second-order unbiased estimator of the mean squared error matrix of the empirical Bayes estimator $\bthh_a^{EB}$ in (\ref{eqn:EB}).
A naive estimator of ${\rm MSEM}(\bthh_a^{EB})$ is the plug-in estimator of (\ref{eqn:MSE}) given by $\G_{1a}(\bPsih^+) + \G_{2a}(\bPsih^+) + \G_{3a}(\bPsih^+) $, but this has a second-order bias, because $E[\G_{1a}(\bPsih^+)]=\G_{1a}(\bPsi)+O(m^{-1})$.
Thus, we need to correct the second-order bias.
Let
\begin{equation}
\G_{4a}(\bPsi)=-\D_a(\bPsi+\D_a)^{-1}{\rm Bias}_{\bPsih}(\bPsi)(\bPsi+\D_a)^{-1}\D_a,
\label{eqn:bias}
\end{equation}
where ${\rm Bias}(\bPsih)$ is the bias of $\bPsih$ given by
$$
{\rm Bias}_{\bPsih}(\bPsi)=
\left\{ \begin{array}{ll} {\rm Bias}_{\bPsih_0}(\bPsi) & {\rm for}\ \bPsih=\bPsih_0,\\
\zero & {\rm for}\ \bPsih=\bPsih_1,
\end{array}\right.
$$
where ${\rm Bias}_{\bPsih_0}(\bPsi)$ is given in (\ref{eqn:Bias}).
Define the estimator ${\rm msem}(\bthh_a^{EB})$ by
\begin{equation}
{\rm msem}(\bthh_a^{EB})=\G_{1a}(\bPsih^+) + \G_{2a}(\bPsih^+) + 2\G_{3a}(\bPsih^+) +\G_{4a}(\bPsih^+).
\label{eqn:msemest}
\end{equation}

\begin{thm}
\label{thm:msemest}
Under the assumption, $E[\G_{1a}(\bPsih^+)+\G_{3a}(\bPsih^+)+\G_{4a}(\bPsih^+)]=\G_{1a}(\bPsi)+O(m^{-3/2})$, and 
$$
E[{\rm msem}(\bthh_a^{EB})]={\rm MSEM}(\bthh_a^{EB})+O(m^{-3/2}),
$$
namely, ${\rm msem}(\bthh_a^{EB})$ is a second-order unbiased estimator of ${\rm MSEM}(\bthh_a^{EB})$.
\end{thm}

{\bf Proof}.\ \ From (2) in Theorem \ref{thm:BLUP}, it is sufficient to show this approximation for $\bPsih$ instead of $\bPsih^+$.
Using the equation in (\ref{eqn:identity}), we can rewrite $G_{1a}(\bPsih)$ as 
\begin{align}
\G_{1a}(\bPsih)=& (\bPsih^{-1}+\D_a^{-1})^{-1}
=\D_a - \D_a(\bPsih+\D_a)^{-1}\D_a
\non\\
=&
\G_{1a}(\bPsi) + \D_a(\bPsi+\D_a)^{-1}(\bPsih-\bPsi)(\bPsi+\D_a)^{-1}\D_a
\label{eqn:BC0}\\
&- \D_a(\bPsi+\D_a)^{-1}(\bPsih-\bPsi)(\bPsi+\D_a)^{-1}(\bPsih-\bPsi)(\bPsi+\D_a)^{-1}\D_a+O_p(m^{-3/2}).\non
\end{align}
We shall evaluate each term in RHS of the above equality.
It is easy to see from (\ref{eqn:Psiap}) that $E[\bPsih-\bPsi]={\rm Bias}(\bPsih)$, which is written as (\ref{eqn:Bias}).
We next evaluate $E[(\bPsih-\bPsi)(\bPsi+\D_a)^{-1}(\bPsih-\bPsi)]$, which is, from (\ref{eqn:Psiap}),  approximated as
\begin{align*}
E[(&\bPsih-\bPsi)(\bPsi+\D_a)^{-1}(\bPsih-\bPsi)]\\
=&{1\over m^2}\sum_{i=1}^m\sum_{j=1}^m E\Big[\Big\{(\y_i-\X_i\bbe)(\y_i-\X_i\bbe)^\top - (\bPsi+\D_i)\Big\}(\bPsi+\D_a)^{-1}
\\
&\times \Big\{(\y_j-\X_j\bbe)(\y_j-\X_j\bbe)^\top - (\bPsi+\D_j)\Big\}\Big] + O(m^{-3/2})
\\
=&
{1\over m^2}\sum_{i=1}^m\sum_{j=1}^m (\bSi+\D_i)^{1/2}
E[( \z_i\z_i^\top-\I) \C ( \z_j\z_j^\top-\I)]  (\bPsi+\D_j)^{1/2} + O(m^{-3/2}),
\end{align*}
for $\C=(\bPsi+\D_i)^{1/2} (\bPsi+\D_a)^{-1} (\bPsi+\D_i)^{1/2}$.
For $i\not= j$, $E[( \z_i\z_i^\top-\I) \C ( \z_j\z_j^\top-\I)] =0$, we have 
$$
\sum_{i=1}^m\sum_{j=1}^m E[( \z_i\z_i^\top-\I) \C ( \z_j\z_j^\top-\I)] 
=\sum_{i=1}^m E[\z_i\z_i^\top\C\z_i\z_i^\top-\C].
$$
Because $E[\z_i\z_i^\top\C\z_i\z_i^\top]=2\C+(\tr\C)\I_k$, it is concluded that 
$$
\D_a(\bPsi+\D_a)^{-1}E[(\bPsih-\bPsi)(\bPsi+\D_a)^{-1}(\bPsih-\bPsi)](\bPsi+\D_a)^{-1}\D_a
=\G_{3a}(\bPsi)+O(m^{-3/2}).
$$

The above arguments imply that a second-order unbiased estimator of $\G_{1a}(\bPsi)$ is $\G_{1a}(\bPsih^+) + \G_{3a}(\bPsih^+) +\G_{4a}(\bPsih^+)$.
The estimators $\G_{2a}(\bPsih^+)$ and $\G_{3a}(\bPsih^+)$ do not have second-order biases, and the results in Theorem \ref{thm:msemest} are established.
\hfill$\Box$

\section{Simulation and Empirical  Studies}
\label{sec:sim}

\subsection{Finite sample performances}

We now investigate finite sample performances of EBLUP in terms of MSEM and the second-order unbiased estimator of MSEM by simulation.

\bigskip
{\bf [1] \ Setup of simulation experiments}.\ \ 
We treat the multivariate Fay-Herriot model (\ref{eqn:MFH}) for $k=2, 3$ and $m=30, 60$ without covariates, namely $\X_i =\I_{k}$.
As a setup of the covariance matrix $\bPsi$ of the random effects, we consider 
$$
\bPsi= \left\{\begin{array}{ll}
\rho\bpsi_2\bpsi_2^\top+(1-\rho){\rm diag}(\bpsi_2\bpsi_2^\top) & {\rm for}\ k=2,\\
\rho\bpsi_3\bpsi_3^\top+(1-\rho){\rm diag}(\bpsi_3\bpsi_3^\top) & {\rm for}\ k=3,
\end{array}\right.
$$
where $\bpsi_2=(\sqrt{1.5}, \sqrt{0.5})^\top$, $\bpsi_3=(\sqrt{1.5}, 1, \sqrt{0.5})^\top$, and ${\rm diag}(\A)$ denotes the diagonal matrix consisting of diagonal elements of matrix $\A$.
Here, $\rho$ is the correlation coefficient, and  we handle the three cases $\rho=0.25, 0.5, 0.75$.
The cases of negative correlations are omitted, because we observe the same results with those of positive ones.

\medskip
Concerning the dispersion matrices $\D_i$ of sampling errors $\bep_i$, we treat the two $\D_i$-patterns: (a) $0.7{\bf I}_k$, $0.6{\bf I}_k$, $0.5{\bf I}_k$, $0.4{\bf I}_k$, $0.3{\bf I}_k$ and (b) $2.0{\bf I}_k$, $0.6{\bf I}_k$, $0.5{\bf I}_k$, $0.4{\bf I}_k$, $0.2{\bf I}_k$.
In the univariate Fay-Herriot model, these cases are treated by Datta, $\et$ (2005).
There are five groups $G_1, \ldots, G_5$ corresponding to these $\D_i$-patterns, 
 and there are six and twelve small areas in each group for $m=30$ and $60$, respectively, where the sampling covariance matrices $\D_i$ are the same for areas within the same group.

\bigskip
{\bf [2]\ Comparison of MSEM}.\ \ 
We begin with obtaining the true mean squared error matrices of the EBLUP $\bthh_a^{EB} = \bthh_a(\bPsih^+)$ by simulation.
Let $\{\y_i^{(r)}, i=1, \ldots, m\}$ be the simulated data in the $r$-th replication for $r=1,\ldots, R$ with $R=50,000$.
Let $\bPsih^{+(r)}$ and $\bth_a^{(r)}$ be the values of $\bPsih^+$ and $\bth_a$ in the $r$-th replication.
Then the simulated value of the true mean squared error matrices is calculated by
$$
{\rm MSEM}(\bthh_a^{EB}) = R^{-1}\sum_{i=1}^R \big\{\bthh_a(\bPsih^{+(r)})-\bth_a^{(r)}\big\}\big\{\bthh_a(\bPsih^{+(r)})-\bth_a^{(r)}\big\}^\top.
$$
As an estimator of $\bPsi$, we here use the simple estimator $\bPsih_0^+$, because  there is little difference between $\bPsih_0^+$ and $\bPsih_1^+$ in simulated values of MSEM under the setup of $\X_i = {\bf I}_{k}$.
Simulated values of the mean squared error matrices, averaged over areas within groups $G_t$, are reported in Tables \ref{table:msem1},  \ref{table:msem2}, and \ref{table:msem3}.
To measure relative improvement of EBLUP, we calculate the percentage relative improvement in the average loss (PRIAL) of $\bthh_a^{EB}$ over $\y_a$, defined by
$$
{\rm PRIAL}(\bthh_a^{EB}, \y_a )
= 100 \times \Big[ 1 - { \tr\{ {\rm MSEM}(\bthh_a^{EB})\} \over \tr\{ {\rm MSEM}(\y_a)\} } \Big].
$$
It is also interesting to compare $\bthh_a^{EB}$ with the EBLUP $\bthh_a^{uEB}$ derived from the univariate Fay-Herriot model.
Thus, we calculate the PRIAL given by
$$
{\rm PRIAL}(\bthh_a^{EB}, \bthh_a^{uEB} )
= 100 \times \Big[ 1 - { \tr\{ {\rm MSEM}(\bthh_a^{EB})\} \over \tr\{ {\rm MSEM}(\bthh_a^{uEB})\} } \Big],
$$
and those values are reported in Tables \ref{table:ir1}, \ref{table:ir2} and \ref{table:ir3}.

\medskip
Table \ref{table:msem1} reports the simulated values of the true MSEM of $\bthh_a^{EB}$ for $k=2$, $\D_i$-patterns (a), $m=30, 60$ and $\rho=0.25, 0.5, 0.75$.
For fixed $m$, the values of MSEM decrease as the correlation $\rho$ in the random effect becomes large.
For fixed $\rho$, the values of MSEM decrease as $m$ becomes large.
Table \ref{table:ir1} reports the values of PRIAL of $\bthh_a^{EB}$ over $\y_a$ and $\bthh_a^{uEB}$ under the same setup as in Table \ref{table:msem1}.  
In all the cases, $\bthh_a^{EB}$ improves on $\y_a$ largely and the improvement rates are larger for larger $\rho$.
In comparison with $\bthh_a^{uEB}$, the univariate EBLUP $\bthh_a^{uEB}$ is slightly better than $\bthh_a^{EB}$ for $\rho=0.25$, but the difference is not significant.
The values of PRIAL of $\bthh_a^{EB}$ over $\bthh_a^{uEB}$ get larger as $\rho$ increases.
In the case of $m=60$, the improvements of $\bthh_a^{EB}$ in light of PRIAL get larger for larger $\rho$.
In the case of $\rho=0.25$, the improvement of $\bthh_a^{EB}$ over $\bthh_a^{uEB}$ is better for $m=60$ than for $m=30$.
This is because the low accuracy in estimation of the covariance matrix $\bPsi$ has more adverse influence on prediction than the benefit from incorporating the small correlation into the estimation.

\medskip 
The comparison of performances between $\D_i$-patterns (a) and (b) is investigated in Tables \ref{table:msem2} and \ref{table:ir2}.
The simulated values of the MSEM of $\bthh_a^{EB}$ in $\D_i$-patterns (a) and (b) are reported in Table \ref{table:msem2}  for $k=2$, $m=30, 60$ and $\rho=0.5$.
As the increment of variance of sampling error in $G_1$, the MSEM in $G_1$ becomes larger, and the other groups have slightly larger MSEM except $G_5$.
The values of PRIAL of $\bthh_a^{EB}$ over $\y_a$ and $\bthh_a^{uEB}$ are given in Table \ref{table:ir2} for $\D_i$-patterns (a) and (b).
under the same setup as in Table \ref{table:msem2}.
As seen from the table, the improvement of $\bthh_a^{EB}$ over $\y_a$ in $G_1$ is larger for $\D_i$-pattern (b) because of the large sampling variance.
However, $\bthh_a^{EB}$ is not better than $\bthh^{uEB}$ in $G_4$ and $G_5$ for $m=30$ and in $G_5$ for $m=60$ in $\D_i$-pattern (b).
This implies that incorporating the information of areas with large sampling variances affects more adversely estimation of areas with small sampling variances in the multivariate model than in the univariate model.

\medskip
Tables \ref{table:msem3} and \ref{table:ir3} report the values of MSEM and PRIAL for $k=3$, $m=30$, $\rho=0.5$ and  $\D_i$-pattern (a).
From Table \ref{table:ir3}, it is revealed that PRIAL of $\bthh_a^{EB}$ over $\y_a$ and $\bthh_a^{uEB}$ are larger for $k=3$ than for $k=2$ in the case of $\rho=0.75$, but smaller in the case of $\rho=0.25$.
When $m$ is fixed as $m=30$, the accuracy in estimation of the covariance matrix $\bPsi$ gets smaller for the larger dimension.
This demonstrates that it is not appropriate to treat the multivariate Fay-Herriot model with a large covariance matrix when $m$ is not large.

\begin{table}[htbp]
\caption{\scriptsize{Simulated values of mean squared error matrices of $\bthh_a^{EB}$ multipled by $100$ for $k=2$, $\D_i$-patterns (a)}}
\label{table:msem1}
\begin{center}
\resizebox{11cm}{!} {
\begin{tabular}{cccc}
&&$m=30$&\\
\hline
 & $\rho=0.25$ & $\rho=0.5$ & $\rho=0.75$ \\
$G_1$ & $
\left[
\begin{array}{rr}
49.8 & 3.8\\
3.8 & 32.6
\end{array}
\right]
$
&
$
\left[
\begin{array}{rr}
48.7 & 8.1\\
8.1 & 30.1
\end{array}
\right]
$
&
$
\left[
\begin{array}{rr}
46.5 & 13.8\\
13.8 & 25.3
\end{array}
\right]
$
\\

$G_2$ & $
\left[
\begin{array}{rr}
44.7 & 3.1\\
3.1 & 30.4
\end{array}
\right]
$
&
$
\left[
\begin{array}{rr}
43.8 & 6.5\\
6.5 & 28.3
\end{array}
\right]
$
&
$
\left[
\begin{array}{rr}
41.4 & 11.6\\
11.6 & 23.7
\end{array}
\right]
$
\\

$G_3$ & $
\left[
\begin{array}{rr}
39.0 & 2.4\\
2.4 & 27.9
\end{array}
\right]
$
&
$
\left[
\begin{array}{rr}
38.0 & 5.3\\
5.3 & 26.3
\end{array}
\right]
$
&
$
\left[
\begin{array}{rr}
36.6 & 9.2\\
9.2 & 21.8
\end{array}
\right]
$
\\

$G_4$ & $
\left[
\begin{array}{rr}
33.1 & 1.7\\
1.7 & 25.3
\end{array}
\right]
$
&
$
\left[
\begin{array}{rr}
32.4 & 3.8\\
3.8 & 23.6
\end{array}
\right]
$
&
$
\left[
\begin{array}{rr}
30.6 & 6.8\\
6.8 & 19.8
\end{array}
\right]
$
\\

$G_5$ & $
\left[
\begin{array}{rr}
26.1 & 1.1\\
1.1 & 21.6
\end{array}
\right]
$
&
$
\left[
\begin{array}{rr}
25.6 & 2.3\\
2.3 & 20.4
\end{array}
\right]
$
&
$
\left[
\begin{array}{rr}
24.2 & 4.6\\
4.6 & 17.4
\end{array}
\right]
$
\\
&&$m=60$&\\
\hline
 & $\rho=0.25$ & $\rho=0.5$ & $\rho=0.75$ \\
$G_1$ & $
\left[
\begin{array}{rr}
49.0 & 4.1\\
4.1 & 30.7
\end{array}
\right]
$
&
$
\left[
\begin{array}{rr}
47.4 & 8.2\\
8.2 & 28.0
\end{array}
\right]
$
&
$
\left[
\begin{array}{rr}
45.2 & 14.0\\
14.0 & 23.6
\end{array}
\right]
$
\\

$G_2$ & $
\left[
\begin{array}{rr}
43.5 & 3.4\\
3.4 & 28.6
\end{array}
\right]
$
&
$
\left[
\begin{array}{rr}
42.5 & 7.0\\
7.0 & 26.5
\end{array}
\right]
$
&
$
\left[
\begin{array}{rr}
40.3 & 11.7\\
11.7 & 22.1
\end{array}
\right]
$
\\

$G_3$ & $
\left[
\begin{array}{rr}
37.9 & 2.6\\
2.6 & 26.0
\end{array}
\right]
$
&
$
\left[
\begin{array}{rr}
37.1 & 5.7\\
5.7 & 24.5
\end{array}
\right]
$
&
$
\left[
\begin{array}{rr}
35.2 & 9.6\\
9.6 & 20.4
\end{array}
\right]
$
\\

$G_4$ & $
\left[
\begin{array}{rr}
31.9 & 1.8\\
1.8 & 23.4
\end{array}
\right]
$
&
$
\left[
\begin{array}{rr}
31.4 & 4.1\\
4.1 & 21.8
\end{array}
\right]
$
&
$
\left[
\begin{array}{rr}
29.8 & 7.3\\
7.3 & 18.5
\end{array}
\right]
$
\\

$G_5$ & $
\left[
\begin{array}{rr}
25.2 & 1.2\\
1.2 & 19.8
\end{array}
\right]
$
&
$
\left[
\begin{array}{rr}
24.8 & 2.7\\
2.7 & 18.7
\end{array}
\right]
$
&
$
\left[
\begin{array}{rr}
23.8 & 5.1\\
5.1 & 16.1
\end{array}
\right]
$
\\
\hline
\end{tabular}
}
\end{center}
\end{table}

\begin{table}[htbp]
\caption{\scriptsize{PRIAL of $\bthh_a^{EB}$ over $\y_a$ and $\bthh_a^{uEB}$ for $k=2$, $\D_i$-patterns (a)}}
\label{table:ir1}
\begin{center}
\resizebox{11cm}{!} {
\begin{tabular}{ccccccc}
&& $\bthh_a^{EB}$ vs $\y_a$ &&& $\bthh_a^{EB}$ vs $\bthh_a^{uEB}$ &\\
\hline
$m=30$ & $\rho=0.25$ & $\rho=0.5$ & $\rho=0.75$ & $\rho=0.25$ & $\rho=0.5$ & $\rho=0.75$\\
$G_1$ & 41.2 & 43.8 & 48.9 & -0.5 & 3.8 & 11.6\\
$G_2$ & 37.2 & 40.1 & 45.7 & 0.0 & 3.5 & 12.3\\
$G_3$ & 33.0 & 35.8 & 41.8 & -0.7 & 3.4 & 11.8\\
$G_4$ & 27.3 & 29.8 & 37.2 & -1.9 & 1.8 & 11.0\\
$G_5$ & 20.8 & 23.5 & 30.4 & -2.5 & 1.1 & 10.0\\
&& $\bthh_a^{EB}$ vs $\y_a$ &&& $\bthh_a^{EB}$ vs $\bthh_a^{uEB}$ &\\
\hline
$m=60$ & $\rho=0.25$ & $\rho=0.5$ & $\rho=0.75$ & $\rho=0.25$ & $\rho=0.5$ & $\rho=0.75$\\
$G_1$ & 43.2 & 45.8 & 51.0 & -0.6 & 4.6 & 13.9\\
$G_2$ & 39.8 & 42.4 & 48.1 & 0.2 & 5.1 & 13.6\\
$G_3$ & 35.6 & 38.7 & 44.2 & 1.3 & 4.6 & 14.3\\
$G_4$ & 30.6 & 33.7 & 39.8 & 0.3 & 3.5 & 13.2\\
$G_5$ & 24.8 & 27.5 & 33.6 & 0.4 & 2.9 & 11.0\\
\hline
\end{tabular}
}
\end{center}
\end{table}

\begin{table}[htbp]
\caption{\scriptsize{Simulated values of mean squared error matrices of $\bthh_a^{EB}$ multipled by $100$ for $k=2$, $\rho=0.5$}}
\label{table:msem2}
\begin{center}
\resizebox{15cm}{!} {
\begin{tabular}{cccccc}
\hline
$m=30$ & Pattern (a) & Pattern (b) & $m=60$ & Pattern (a) & Pattern (b) \\
$G_1$ & $
\left[
\begin{array}{rr}
48.7 & 8.1\\
8.1 & 30.1
\end{array}
\right]
$
&
$
\left[
\begin{array}{rr}
89.9 & 19.7\\
19.7 & 42.9
\end{array}
\right]
$
&$G_1$ & $
\left[
\begin{array}{rr}
47.4 & 8.2\\
8.2 & 28.0
\end{array}
\right]
$
&
$
\left[
\begin{array}{rrr}
86.8 & 20.1\\
20.1 & 40.0
\end{array}
\right]
$
\\

$G_2$ & $
\left[
\begin{array}{rr}
43.8 & 6.5\\
6.5 & 28.3
\end{array}
\right]
$
&
$
\left[
\begin{array}{rr}
44.5 & 6.0\\
6.0 & 30.2
\end{array}
\right]
$
&$G_2$ & $
\left[
\begin{array}{rr}
42.5 & 7.0\\
7.0 & 26.5
\end{array}
\right]
$
&
$
\left[
\begin{array}{rr}
42.9 & 6.5\\
6.5 & 27.8
\end{array}
\right]
$
\\

$G_3$ & $
\left[
\begin{array}{rr}
38.0 & 5.3\\
5.3 & 26.3
\end{array}
\right]
$
&
$
\left[
\begin{array}{rrr}
39.3 & 4.7\\
4.7 & 28.3
\end{array}
\right]
$
&$G_3$ & $
\left[
\begin{array}{rr}
37.1 & 5.7\\
5.7 & 24.5
\end{array}
\right]
$
&
$
\left[
\begin{array}{rrr}
37.8 & 5.0\\
5.0 & 25.8
\end{array}
\right]
$
\\

$G_4$ & $
\left[
\begin{array}{rr}
32.4 & 3.8\\
3.8 & 23.6
\end{array}
\right]
$
&
$
\left[
\begin{array}{rrr}
33.4 & 3.2\\
3.2 & 25.9
\end{array}
\right]
$
&$G_4$ & $
\left[
\begin{array}{rr}
31.4 & 4.1\\
4.1 & 21.8
\end{array}
\right]
$
&
$
\left[
\begin{array}{rrr}
32.0 & 3.6\\
3.6 & 23.8
\end{array}
\right]
$
\\

$G_5$ & $
\left[
\begin{array}{rr}
25.6 & 2.3\\
2.3 & 20.4
\end{array}
\right]
$
&
$
\left[
\begin{array}{rrr}
19.1 & 0.1\\
0.1 & 18.8
\end{array}
\right]
$
&$G_5$ & $
\left[
\begin{array}{rr}
24.8 & 2.7\\
2.7 & 18.7
\end{array}
\right]
$
&
$
\left[
\begin{array}{rrr}
18.1 & 0.6\\
0.6 & 16.4
\end{array}
\right]
$

\\
\hline
\end{tabular}
}
\end{center}
\end{table}

\begin{table}[htbp]
\caption{\scriptsize{PRIAL of $\bthh_a^{EB}$ over $\y_a$ and $\bthh_a^{uEB}$ for $k=2$, $m=30, 60$, $\rho=0.5$, $\D_i$-patterns (a), (b)}}
\label{table:ir2}
\begin{center}
\resizebox{11cm}{!} {
\begin{tabular}{ccccc}
& $\bthh_a^{EB}$ vs $\y_a$ && $\bthh_a^{EB}$ vs $\bthh_a^{uEB}$ &\\
\hline
$m=30$ & Pattern (a) & Pattern (b) & Pattern (a) & Pattern (b)\\
$G_1$ & 43.8 & 66.4 & 3.8 & 2.1\\
$G_2$ & 40.1 & 37.0 & 3.5 & 0.8\\
$G_3$ & 35.8 & 32.1 & 3.4 & 1.0\\
$G_4$ & 29.8 & 26.2 & 1.8 & -0.2\\
$G_5$ & 23.5 & 4.2 & 1.1 & -8.5\\
& $\bthh_a^{EB}$ vs $\y_a$ && $\bthh_a^{EB}$ vs $\bthh_a^{uEB}$ &\\
\hline
$m=60$ & Pattern (a) & Pattern (b) & Pattern (a) & Pattern (b)\\
$G_1$ & 45.8 & 68.5 & 4.6 & 3.1\\
$G_2$ & 42.4 & 40.7 & 5.1 & 3.1\\
$G_3$ & 38.7 & 36.2 & 4.6 & 2.7\\
$G_4$ & 33.7 & 30.6 & 3.5 & 1.3\\
$G_5$ & 27.5 & 13.9 & 2.9 & -2.7\\
\hline
\end{tabular}
}
\end{center}
\end{table}

\begin{table}[htbp]
\caption{\scriptsize{Simulated values of mean squared error matrices of $\bthh_a^{EB}$ multiplied by $100$ for $k=3$, $m=30$, $\D_i$-patterns (a)}}
\label{table:msem3}
\begin{center}
\resizebox{11cm}{!} {
\begin{tabular}{cccc}
&&$m=30$&\\
\hline
 & $\rho=0.25$ & $\rho=0.5$ & $\rho=0.75$ \\
G1 & $
\left[
\begin{array}{rrr}
50.0 & 3.4 & 3.5\\
3.4 & 44.3 & 3.4\\
3.5 & 3.4 & 33.3
\end{array}
\right]
$
&
$
\left[
\begin{array}{rrr}
48.0 & 7.0 & 6.4\\
7.0 & 41.1 & 6.5\\
6.4 & 6.5 & 29.5
\end{array}
\right]
$
&
$
\left[
\begin{array}{rrr}
42.0 & 12.3 & 10.0\\
12.3 & 34.8 & 9.4\\
10.0 & 9.4 & 23.1
\end{array}
\right]
$
\\

G2 & $
\left[
\begin{array}{rrr}
45.2 & 2.6 & 2.8\\
2.6 & 39.8 & 2.9\\
2.8 & 2.9 & 31.2
\end{array}
\right]
$
&
$
\left[
\begin{array}{rrr}
42.8 & 5.8 & 5.4\\
5.8 & 37.7 & 5.5\\
5.4 & 5.5 & 28.2
\end{array}
\right]
$
&
$
\left[
\begin{array}{rrr}
38.2 & 10.0 & 8.3\\
10.0 & 31.8 & 8.0\\
8.3 & 8.0 & 21.7
\end{array}
\right]
$
\\

G3 & $
\left[
\begin{array}{rrr}
40.0 & 2.0 & 1.9\\
1.9 & 36.1 & 2.1\\
1.9 & 2.1 & 29.0
\end{array}
\right]
$
&
$
\left[
\begin{array}{rrr}
37.5 & 4.1 & 3.9\\
4.1 & 33.9 & 4.1\\
3.9 & 4.1 & 25.8
\end{array}
\right]
$
&
$
\left[
\begin{array}{rrr}
33.5 & 7.7 & 7.0\\
7.7 & 28.8 & 6.4\\
7.0 & 6.4 & 20.5
\end{array}
\right]
$
\\

G4 & $
\left[
\begin{array}{rrr}
33.4 & 1.3 & 1.5\\
1.3 & 31.0 & 1.6\\
1.5 & 1.6 & 26.0
\end{array}
\right]
$
&
$
\left[
\begin{array}{rrr}
32.1 & 2.7 & 2.9\\
2.7 & 29.2 & 3.0\\
2.9 & 3.0 & 20.7
\end{array}
\right]
$
&
$
\left[
\begin{array}{rrr}
29.2 & 5.6 & 5.1\\
5.6 & 25.2 & 5.1\\
5.1 & 5.1 & 18.4
\end{array}
\right]
$
\\

G5 & $
\left[
\begin{array}{rrr}
26.3 & 0.7 & 0.7\\
0.7 & 25.4 & 1.0\\
0.7 & 1.0 & 22.7
\end{array}
\right]
$
&
$
\left[
\begin{array}{rrr}
25.8 & 1.6 & 1.5\\
1.6 & 24.1 & 1.8\\
1.5 & 1.8 & 20.7
\end{array}
\right]
$
&
$
\left[
\begin{array}{rrr}
23.4 & 3.1 & 3.2\\
3.1 & 21.0 & 3.4\\
3.2 & 3.4 & 16.5
\end{array}
\right]
$
\\
\hline
\end{tabular}
}
\end{center}
\end{table}

\begin{table}[htbp]
\caption{\scriptsize{PRIAL of $\bthh_a^{EB}$ over $\y_a$ and $\bthh_a^{uEB}$ for $k=2,3$, $m=30$, $\D_i$-patterns (a)}}
\label{table:ir3}
\begin{center}
\resizebox{11cm}{!} {
\begin{tabular}{ccccccc}
&& $\bthh_a^{EB}$ vs $\y_a$ &&& $\bthh_a^{EB}$ vs $\bthh_a^{uEB}$ &\\
\hline
$k=2$ & $\rho=0.25$ & $\rho=0.5$ & $\rho=0.75$ & $\rho=0.25$ & $\rho=0.5$ & $\rho=0.75$\\
$G_1$ & 41.2 & 43.8 & 48.9 & -0.5 & 3.8 & 11.6\\
$G_2$ & 37.2 & 40.1 & 45.7 & 0.0 & 3.5 & 12.3\\
$G_3$ & 33.0 & 35.8 & 41.8 & -0.7 & 3.4 & 11.8\\
$G_4$ & 27.3 & 29.8 & 37.2 & -1.9 & 1.8 & 11.0\\
$G_5$ & 20.8 & 23.5 & 30.4 & -2.5 & 1.1 & 10.0\\
&& $\bthh_a^{EB}$ vs $\y_a$ &&& $\bthh_a^{EB}$ vs $\bthh_a^{uEB}$ &\\
\hline
$k=3$ & $\rho=0.25$ & $\rho=0.5$ & $\rho=0.75$ & $\rho=0.25$ & $\rho=0.5$ & $\rho=0.75$\\
G1 & 39.5 & 43.6 & 52.4 & -1.9 & 5.9 & 20.3\\
G2 & 35.2 & 40.0 & 48.9 & -2.6 & 4.8 & 19.2\\
G3 & 30.3 & 35.1 & 44.6 & -3.9 & 3.2 & 18.2\\
G4 & 25.2 & 29.7 & 39.6 & -4.4 & 1.9 & 15.7\\
G5 & 17.6 & 21.6 & 32.0 & -5.7 & -0.3 & 12.8\\
\hline
\end{tabular}
}
\end{center}
\end{table}

\bigskip
{\bf [3]\ MSEM approximation and its estimator}.\ \ 
We next investigate the performance of the second-order approximation of MSEM of EBLUP $\bthh_a^{EB}$ given in Theorem \ref{thm:MSE} and the second-order unbiased estimator ${\rm msem}(\bthh_a^{EB})$ of MSEN given in Theorem \ref{thm:msemest}.
The values of the second-order approximation of MSEM are given in Table \ref{table:soa1} for $k=2$, $m=3$ and $\D_i$-pattern (a).
Comparing the values in Table \ref{table:soa1} with the corresponding true values of the MSEM  in Table \ref{table:msem1}, we can see that the second-order approximation can approximate the true MSEM precisely for every $G_t$ and $\rho$.

\medskip
Concerning the performance of the second-order unbiased estimator ${\rm msem}(\bthh_a^{EB})$ given in (\ref{eqn:msemest}), we compute the simulated values of relative bias of the estimator ${\rm msem}(\bthh_a^{EB})$, averaged over areas within groups $G_t$.
Those values are reported in Table \ref{table:rb1} for $k=2$, $m=30, 60$ and $\D_i$-pattern (a).
It is revealed from Table \ref{table:rb1} that the relative bias gets larger for larger $\rho$.
Also, the values of the relative bias are smaller for $m=60$ than for $m=30$, namely, the relative bias gets small as $m$ increases.

\begin{table}[htbp]
\caption{\scriptsize{Second order approximations of mean squared error matrices of $\bthh_a^{EB}$ multiplied by $100$ for $k=2$, $\D_i$-patterns (a)}}
\label{table:soa1}
\begin{center}
\resizebox{11cm}{!} {
\begin{tabular}{cccc}
&&$m=30$&\\
\hline
 & $\rho=0.25$ & $\rho=0.5$ & $\rho=0.75$ \\
$G_1$ & $
\left[
\begin{array}{rr}
49.8 & 3.7\\
3.7 & 32.6
\end{array}
\right]
$
&
$
\left[
\begin{array}{rr}
48.6 & 7.9\\
7.9 & 30.3
\end{array}
\right]
$
&
$
\left[
\begin{array}{rr}
46.2 & 13.2\\
13.2 & 25.9
\end{array}
\right]
$
\\

$G_2$ & $
\left[
\begin{array}{rr}
44.6 & 3.1\\
3.1 & 30.4
\end{array}
\right]
$
&
$
\left[
\begin{array}{rr}
43.6 & 6.6\\
6.6 & 28.4
\end{array}
\right]
$
&
$
\left[
\begin{array}{rr}
41.5 & 11.1\\
11.1 & 24.4
\end{array}
\right]
$
\\

$G_3$ & $
\left[
\begin{array}{rr}
38.9 & 2.4\\
2.4 & 27.8
\end{array}
\right]
$
&
$
\left[
\begin{array}{rr}
38.1 & 5.2\\
5.2 & 26.1
\end{array}
\right]
$
&
$
\left[
\begin{array}{rr}
36.3 & 8.9\\
8.9 & 22.6
\end{array}
\right]
$
\\

$G_4$ & $
\left[
\begin{array}{rr}
32.6 & 1.7\\
1.7 & 24.7
\end{array}
\right]
$
&
$
\left[
\begin{array}{rr}
32.0 & 3.8\\
3.8 & 23.3
\end{array}
\right]
$
&
$
\left[
\begin{array}{rr}
30.6 & 6.6\\
6.6 & 20.5
\end{array}
\right]
$
\\

$G_5$ & $
\left[
\begin{array}{rr}
25.7 & 1.1\\
1.1 & 20.7
\end{array}
\right]
$
&
$
\left[
\begin{array}{rr}
25.3 & 2.4\\
2.4 & 20.0
\end{array}
\right]
$
&
$
\left[
\begin{array}{rr}
24.4 & 4.3\\
4.3 & 17.8
\end{array}
\right]
$
\\
\hline
\end{tabular}
}
\end{center}
\end{table}

\begin{table}[htbp]
\caption{\scriptsize{Simulated values of percentage average relative bias of mean squared error matrices of $\bthh_a^{EB}$ multiplied by $100$ for $k=2$, $m=30, 60$, $\D_i$-pattern (a)}}
\label{table:rb1}
\begin{center}
\resizebox{11cm}{!} {
\begin{tabular}{cccc}
&&Pattern (a)&\\
\hline
$m=30$ & $\rho=0.25$ & $\rho=0.5$ & $\rho=0.75$ \\
$G_1$ & $
\left[
\begin{array}{rr}
-0.3 & -2.6\\
-2.6 & 1.1
\end{array}
\right]
$
&
$
\left[
\begin{array}{rr}
-0.9 & -4.1\\
-4.1 & 2.9
\end{array}
\right]
$
&
$
\left[
\begin{array}{rr}
0.6 & -9.5\\
-9.5 & 10.1
\end{array}
\right]
$
\\

$G_2$ & $
\left[
\begin{array}{rr}
0.6 & 3.1\\
3.1 & 0.9
\end{array}
\right]
$
&
$
\left[
\begin{array}{rr}
0.3 & -3.5\\
-3.5 & 2.7
\end{array}
\right]
$
&
$
\left[
\begin{array}{rr}
1.1 & -10.4\\
-10.4 & 13.1
\end{array}
\right]
$
\\

$G_3$ & $
\left[
\begin{array}{rr}
-0.6 & -5.8\\
-5.8 & 1.2
\end{array}
\right]
$
&
$
\left[
\begin{array}{rr}
1.3 & -7.8\\
-7.8 & 4.6
\end{array}
\right]
$
&
$
\left[
\begin{array}{rr}
1.2 & -16.7\\
-16.7 & 13.6
\end{array}
\right]
$
\\
$G_4$ & $
\left[
\begin{array}{rr}
-0.4 & -4.6\\
-4.6 & 2.9
\end{array}
\right]
$
&
$
\left[
\begin{array}{rr}
0.4 & -10.8\\
-10.8 & 4.7
\end{array}
\right]
$
&
$
\left[
\begin{array}{rr}
1.2 & -23.4\\
-23.4 & 17.8
\end{array}
\right]
$
\\

$G_5$ & $
\left[
\begin{array}{rr}
0.3 & -24.4\\
-24.4 & 2.2
\end{array}
\right]
$
&
$
\left[
\begin{array}{rr}
0.6 & -26.1\\
-26.1 & 7.7
\end{array}
\right]
$
&
$
\left[
\begin{array}{rr}
3.4 & -42.3\\
-42.3 & 23.1
\end{array}
\right]
$
\\
&&Pattern (a)&\\
\hline
$m=60$ & $\rho=0.25$ & $\rho=0.5$ & $\rho=0.75$ \\
$G_1$ & $
\left[
\begin{array}{rr}
-0.1 & -2.3\\
-2.3 & -0.5
\end{array}
\right]
$
&
$
\left[
\begin{array}{rr}
0.2 & -0.2\\
-0.2 & -0.5
\end{array}
\right]
$
&
$
\left[
\begin{array}{rr}
0.4 & -0.7\\
-0.7 & 1.7
\end{array}
\right]
$
\\

$G_2$ & $
\left[
\begin{array}{rr}
0.7 & -3.4\\
-3.4 & -0.2
\end{array}
\right]
$
&
$
\left[
\begin{array}{rr}
0.4 & -0.8\\
-0.8 & -0.5
\end{array}
\right]
$
&
$
\left[
\begin{array}{rr}
0.8 & -0.0\\
-0.0 & 2.1
\end{array}
\right]
$
\\

$G_3$ & $
\left[
\begin{array}{rr}
0.2 & -5.1\\
-5.1 & -0.2
\end{array}
\right]
$
&
$
\left[
\begin{array}{rr}
0.1 & -1.5\\
-1.5 & 0.1
\end{array}
\right]
$
&
$
\left[
\begin{array}{rr}
0.3 & -1.4\\
-1.4 & 2.9
\end{array}
\right]
$
\\
$G_4$ & $
\left[
\begin{array}{rr}
-0.1 & -5.1\\
-5.1 & -0.4
\end{array}
\right]
$
&
$
\left[
\begin{array}{rr}
-0.4 & -2.4\\
-2.4 & 0.3
\end{array}
\right]
$
&
$
\left[
\begin{array}{rr}
1.4 & -2.9\\
-2.9 & 3.6
\end{array}
\right]
$
\\

$G_5$ & $
\left[
\begin{array}{rr}
0.2 & -3.3\\
-3.3 & -0.2
\end{array}
\right]
$
&
$
\left[
\begin{array}{rr}
0.3 & -5.9\\
-5.9 & -0.1
\end{array}
\right]
$
&
$
\left[
\begin{array}{rr}
0.6 & -8.1\\
-8.1 & 5.1
\end{array}
\right]
$
\\
\hline
\end{tabular}
}
\end{center}
\end{table}

\subsection{Illustrative example}
\label{sec:exm}

This example, primarily for illustration, uses the multivariate Fay-Herriot model (\ref{eqn:MFH}) and data from the 2016 Survey of Family Income and Expenditure in Japan, which is based on two or more person households (excluding agricultural, forestry and fisheries households).
The target domains are the 47 Japanese prefectural capitals.
The 47 prefectures are divided into 10 regions: Hokkaido, Tohoku, Kanto, Hokuriku, Tokai, Kinki, Chugoku, Shikoku, Kyushu and Okinawa.
Each region consists of several prefectures except Hokkaido and Okinawa, which consist of one prefecture.

\medskip
In this study, as observations $(y_{i1}, y_{i2})^\top$, we use the reported data of the yearly averaged monthly spendings on ^^ Education' and ^^  Cultural-amusement' per worker's household, scaled by 1,000 Yen, at each capital city of 47 prefectures.
In addition, we use the data in the 2014 National Survey of Family Income and Expenditure.
The average spending data in this survey are more reliable than the Survey of Family Income and Expenditure since the sample sizes are much larger.
However, this survey is conducted only once in every five years.
As auxiliary variables, we use the data of the average spendings on ^^ Education' and ^^  Cultural-amusement', which is denoted by ${\rm EDU}_{i}$ and ${\rm CUL}_{i}$, respectively.
Then the regressor in the model (\ref{eqn:MFH}) is 
$$
\X_i = \begin{pmatrix} 1 & {\rm EDU}_{i} & 0 & 0 \\ 0 & 0 & 1 & {\rm CUL}_{i}\end{pmatrix}.
$$
Then we apply the multivariate Fay-Herriot model (\ref{eqn:MFH}), where sampling covariance matrices $\D_i$ of the $i$-th region for $i=1,\ldots,10$ are calculated based on data of yearly averaged monthly spendings on ^^ Education' and ^^  Cultural-amusement' in the past ten years (2006-2015), where $\D_i$ is given as the average of the sampling covariance matrices of prefectures within the $i$-th region.
That is,  the sampling covariance matrix $\D_i$ are the same for prefectures within the same region.

\medskip
The estimates of the covariance matrix $\bPsi$ and the correlation coefficient $\rho$ is 
$$
\bPsih=\begin{pmatrix}
8.5 & 3.0\\
3.0 & 10.2
\end{pmatrix}
\quad {\rm and}\quad {\hat \rho}=0.32.
$$
The estimates of the regression coefficients and the $p$-values for testing $H_0: \be_k=0$ for $k=1,\ldots,4$ are given in Table \ref{tab:reg}.
All the estimates are significant.
\begin{table}[htbp]
\caption{\scriptsize{Estimates of regression coefficients and $p$-values}}
\label{tab:reg}
\begin{center}
\resizebox{11cm}{!} {
\begin{tabular}{ccccc}
\hline
varables & Constant(EDU) & EDU & Constant(CUL) & CUL \\
\hline
$\bbe$ & $4.47$ & $0.82$ & $12.12$ & 0.65 \\
$p$-value & $0.007$ & $0.000$ & $0.002$ & $0.000$\\
\hline
\end{tabular}
}
\end{center}
\end{table}

\medskip
The values of EBLUP and direct estimate of spendings on ^^ Education' and ^^  Cultural-amusement' are reported in Table \ref{tab:eblup}.
We only pick up the three prefectures from three different regions: Tokyo prefecture from the Kanto region, Osaka prefecture from the Kinki region and Fukushima prefecture from the Tohoku region, whose sampling covariance matrices are
$$
\begin{pmatrix}
1.1 & 0.3\\
0.3 & 3.0
\end{pmatrix},\ 
\begin{pmatrix}
1.1 & -0.2\\
-0.2 & 3.9
\end{pmatrix}
\quad
{\rm and}
\quad
\begin{pmatrix}
4.7 & 3.5\\
3.5 & 4.9
\end{pmatrix},
$$
respectively.
It is seen that as the sampling variances become larger, the direct estimates are more shrunken by the EBLUP in the sense of (direct estimate - EBLUP)/(direct estimate).

\begin{table}[htbp]
\caption{\scriptsize{EBLUP and direct estimates}}
\label{tab:eblup}
\begin{center}
\resizebox{9cm}{!} {
\begin{tabular}{cccc}
\hline
 & Tokyo & Osaka & Fukushima \\
\hline
direct estimator (EDU) & 32.5 & 19.0 & 13.3 \\
EBLUP (EDU) & 31.8 & 19.4 & 12.6\\
\hline
direct estimator (CUL) & 41.8 & 24.9 & 29.9 \\
EBLUP (CUL) & 40.6 & 26.1 & 29.0\\
\hline
\end{tabular}
}
\end{center}
\end{table}

\medskip
The uncertainty of EBLUP is provided by the second-order unbiased estimator of MSEM of EBLUP.
Table \ref{table:rmse} reports the estimates of MSEM averaged over prefectures within each region for 10 regions.
We also calculate the percentage relative improvement in the average loss estimate (PRIAL estimate) of $\bthh_a^{EB}$ over $\y_a$ and $\bthh_a^{uEB}$. 
Table \ref{table:rir} reports the average of those values over each region for spendings on education and cultural-amusement.
It is revealed from Table \ref{table:rir} that the multivariate EBLUP improves on  the direct estimates significantly and that the multivariate EBLUP is slightly better than the univariate EBLUP for most regions except Okinawa, which has a smaller sampling covariance matrix.

\begin{table}[htbp]
\caption{\scriptsize{Estimates of the mean squared error matrices of $\bthh_a^{EB}$}}
\label{table:rmse}
\begin{center}
\resizebox{13cm}{!} {
\begin{tabular}{ccccc}
\hline
Hokkaido & Tohoku & Kanto & Hokuriku & Tokai \\
$
\left[
\begin{array}{rr}
0.5 & 0.7\\
0.7 & 3.8
\end{array}
\right]
$
&
$
\left[
\begin{array}{rr}
3.2 & 2.2\\
2.2 & 3.4
\end{array}
\right]
$
&
$
\left[
\begin{array}{rr}
1.0 & 0.3\\
0.3 & 2.5
\end{array}
\right]
$
&
$
\left[
\begin{array}{rr}
1.0 & 0.6\\
0.6 & 4.7
\end{array}
\right]
$
&
$
\left[
\begin{array}{rr}
1.4 & 0.6\\
0.6 & 1.8
\end{array}
\right]
$
\\
Kinki & Chugoku & Shikoku & Kyushu & Okinawa \\
$
\left[
\begin{array}{rr}
1.0 & -0.0\\
-0.0 & 2.8
\end{array}
\right]
$
&
$
\left[
\begin{array}{rr}
1.5 & 0.3\\
0.3 & 2.6
\end{array}
\right]
$
&
$
\left[
\begin{array}{rr}
4.2 & 0.9\\
0.9 & 3.5
\end{array}
\right]
$
&
$
\left[
\begin{array}{rr}
1.0 & 0.7\\
0.7 & 1.8
\end{array}
\right]
$
&
$
\left[
\begin{array}{rr}
3.0 & 0.8\\
0.8 & 1.7
\end{array}
\right]
$
\\
\hline
\end{tabular}
}
\end{center}
\end{table}

\begin{table}[htbp]
\caption{\scriptsize{PRIAL estimates of $\bthh_a^{EB}$ over $\y_a$ and $\bthh_a^{uEB}$}}
\label{table:rir}
\begin{center}
\resizebox{17cm}{!} {
\begin{tabular}{cccccccccc}
$\bthh_a^{EB}$ vs $\y_a$ &&&&&&&&&\\
\hline
Hokkaido & Tohoku & Kanto & Hokuriku & Tokai & Kinki & Chugoku & Shikoku & Kyushu & Okinawa\\
82.4 & 84.7 & 80.9 & 84.1 & 80.5 & 81.5 & 81.6 & 85.7 & 80.1 & 81.0\\
\hline
\\
$\bthh_a^{EB}$ vs $\bthh_a^{uEB}$ &&&&&&&&&\\
\hline
Hokkaido & Tohoku & Kanto & Hokuriku & Tokai & Kinki & Chugoku & Shikoku & Kyushu & Okinawa\\
4.1 & 7.1 & 1.9 & 5.3 & 1.1 & 4.5 & 2.4 & 3.9 & 1.3 & -33.9\\
\hline
\end{tabular}
}
\end{center}
\end{table}

\section*{Acknowledgments}
Research of the second author was supported in part by Grant-in-Aid for Scientific Research  (15H01943 and 26330036) from Japan Society for the Promotion of Science.

\end{document}